\documentclass[11pt,twoside]{amsart} \textwidth 12 cm
 \usepackage{graphicx} 

\numberwithin{equation}{section}

\begin{document}

\title{How many points in a point cloud is sufficient for accurate estimation of the curvature }
\author{R. Mirzaie}

\begin{abstract}
We introduce an estimator for the  curvature of curves and
surfaces by using   finite sample points  drawn from sampling a
probability distribution that has support on  the curve or
surface. First we give an algorithm for estimation of
the curvature in a given point of  a curve. Then, we extend it to
estimate the Gaussian curvature of  the surfaces. In the proposed
algorithms, we use a relation between the number of selected
points in the point cloud and the probability that a given point has a sufficient number of
nearby points. This relation allows us  to control the required number of points in the  point cloud.


\end{abstract}
\thanks{{\scriptsize R. Mirzaie, Department of Pure Mathematics,
Faculty of Science, Imam Khomeini International University,
Qazvin, Iran \flushleft  r.mirzaei@sci.ikiu.ac.ir{\scriptsize }}}
\maketitle

\pagestyle{myheadings}

\markboth{\rightline {\scriptsize  R. Mirzaie}}
         {\leftline{Estimation of the curvature from random samples }}
\section{ Introduction}
Recently, there has been significant progress in the manifold
learning and has been attractive for many researches, because of
the useful applications in  topological data analysis, geometric
measure theory, integral geometry, image processing arguments and
biomathematics. But, a main problem is still open: What is the
best algorithm for estimating the curvature of a manifold using a
given  point cloud.
 There are many  algorithms suited for Gaussian curvature and
mean curvature estimation. Almost all of the algorithms rely on
the estimating of intrinsic geometric properties of a surface
from a polygonal mesh obtained from  data. For example, a widely
used approximation of the curvature of a surface is based on the
fact that a surface around a point can be represented as a graph
of a bivariate function, and the curvature can be approximated to
the curvature of a hyperbola passing from the neighbouring points
of the given point in the selected point cloud( see [12]). The
authors of [1] use principal component analysis to estimate the
Ricci curvature of a manifold at each point  of a point cloud set,
 selected randomly from the points of a Riemannian submanifold
of the euclidean space. Estimator of [9], for the scalar
curvature, is based on the fact that the scalar curvature at a
point  $x$ of a Riemannian manifold characterizes the growth rate
of the volume of a geodesic ball $B(x, r)$  as $r$ decreases.
Thus, estimation of  the scalar curvature reduces to estimation
of the volumes. We refer to [ 3, 4, 5, 6, 8, 9, 14] for many methods to
estimate Gaussian curvature from point clouds that are sampled
from surfaces. Also, we refer the reader to [2, 10, 12, 13] for
surveys. The authors of  [7],  locally fit a set of circles through
the neighbour vertices of the triangulations, as a main tool for
estimations.
 In [2] a suitable curve
 is constructed at each vertex of the surface triangulation and
then the local differential structure is derived.
 In the present article,  we give an algorithm to estimate 
 the curvature of a curve by using the osculating
circle at a given point. Then, we use it to estimate the Gussian
curvature of surfaces.\\
Using osculating circle is a traditional method for estimation of the curve of plane curves, but using it for estimation of the principal curvatures
of the surfaces from random point clouds is new. Another subject which is new in the present article is the following fact:\\

Using random point clouds to estimate   geometric properties with global characteristics works very well. But,
when we use random point clouds to estimate  geometric properties with local characteristics (for example curvature at a point), there may be significant estimation errors. Let $a$ be an arbitrary point on a surface $M$ and $B=\{x_{1},...,x_{m}\}$ be a  subset  of points on $M$. To estimate
the curvature at the point $a$,  the points of $B$ that are close to $a$ play essential role. But, we are uncertain if  there are enough points in $B$
close to $a$. For example, in the following picture we have chosen   points randomly on the surface. The selected points  are not suitable
for estimation of the curvature at  $b$. But, they are suitable for the curvature estimation at $a$\\

\noindent\includegraphics[height=8cm]{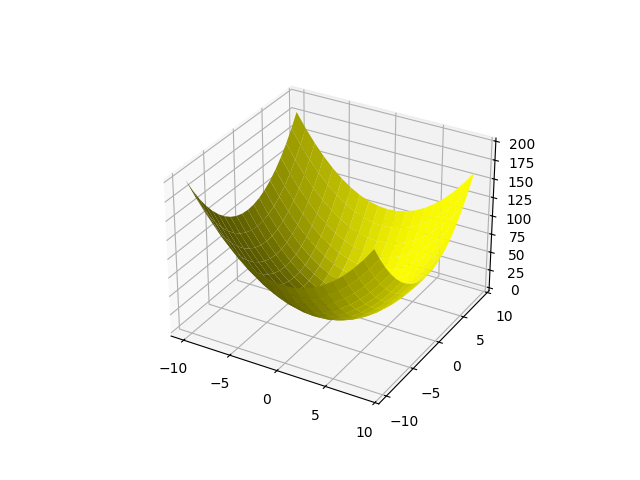}\\

\noindent\includegraphics[height=8cm]{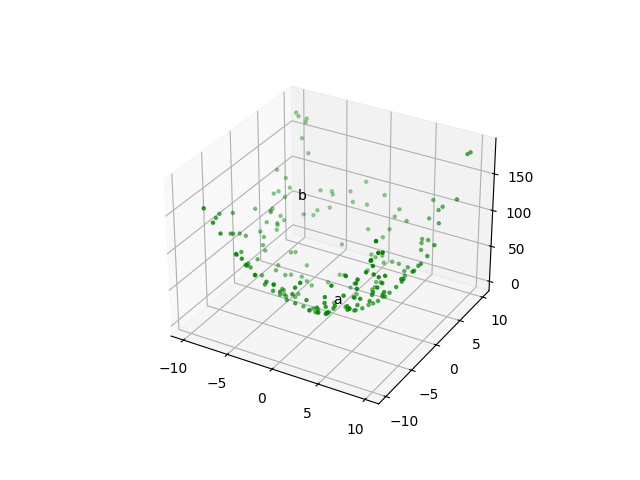}\\
 In practice, since we usually  have no information
about the underlying manifold, the easiest way to solve  the problem is to choose a large number of points. But, the other problem is that the large number of points 
makes the  computation process  and runtime  of the computer program long. 
To overcome these two problems, we use  Theorem 2.7 and Theorem 2.14, to provide a method for finding   
  the minimum number $m$ for  the   points in the point cloud.  $m$   is suitable  with a probability $p$, where $p$   can be chosen arbitrary  close to one.\\\\

\section{Results}
In what follows, each curve and surface is considered to be differentiable  with the finite length and area imbedded in $R^{2}$ and $R^{3}$, respectively. Since we intend to estimate their curvature at a generic point,  we suppose that they have not boundary points. \\
We suppose that each sample set $B=\{x_{1},...,x_{m}\}$ of the points on the curve or surface is i.i.d. (independent and identically distributed) uniformly distributed ( uniform with respect to the Riemannian volume measure induced from the ambient space).\\\\ 
{\bf Remark 2.1.}  Let $\gamma$ be a differentiable plan curve
and let $r=r(a,b,c)$ be  radius of the circle generated by three
nonlinear points $a, b, c$ belonging to $\gamma$.  The curvature
of $\gamma$ at the point $a$ which we denote by $\kappa(a)$ is
defined by
 to $\kappa(a)^{-1}=lim_{_{(b,c) \to (a,a)}} r(a,b,c)$. A circle with radius $\kappa(a)^{-1}$tangent to $\gamma$ at
 the point $a$ is called the osculating circle at $a$.
  Thus, for the points $b, c$ sufficiently close to $a$, $\frac{1}{r}$ is
   an estimation for $\kappa(a)$. If $max\{|b-a|, |c-a|\} < \epsilon$,
   then we say that $\frac{1}{r}$ is a $\epsilon$-triangle estimation of $\kappa(a)$.
   Let $d_{1}, d_{2}, d_{3}$ be the length of the sides of the triangle $abc$ and $D=\frac{d_{1}+d_{2}+d_{3}}{2}$. Then, 
   \[r=\frac{d_{1}d_{2}d_{3}}{4\sqrt{D(D-d_{1})(D-d_{2})(D-d_{3})}}.\]
   Thus, $\epsilon$-triangle estimation of $\kappa(a)$ is 
   \[ \kappa(a)=\frac{4\sqrt{D(D-d_{1})(D-d_{2})(D-d_{3})}}{d_{1}d_{2}d_{3}}.\]
In what follows, we denote by
$d(x,y)$ the distance between $x, y$ in $R^{n}$.\\
If $S$ is a connected curve or a surface, then the distance between $x,y$ in $S$ is denoted by $d_{S}(x,y)$, which is equal to the infimum of the
length of the curves in $S$ joining $x$ to $y$.\\
We denote by $\gamma$ and $M$ a differentiable plan curve with the
finite length and  a differentiable surface with the finite area,
respectively.  $B(x,r)$ is a  ball with the center $x$ and of
radius $r>0$ ( in $R^{n}$). $B(r)$ will denote a ball of radius
$r$ ( non-specified center). If $M$ is a surface in $R^{3}$ and
$x \in M$, we denote by $B_{M}(x,r)$ the geodesic ball in $M$
with the center $x$ and the
 radius $r$, that is \[B_{M}(x,r)=\{ y\in M:  d_{M}(x,y) \leq r\}.\]
We use the symbol $vol(M)$ to denote the volume(area) of $M$.\\\\

{\bf Remark 2.2.} $vol(B_{M}(x,r))$ is a continuous  function of
 $r$. Thus, for  each positive number  $\epsilon$ ( sufficiently small ),
there is a number $r=r(x, \epsilon)$ such that $vol
(B_{M}(x,r))=\pi \epsilon^{2}$ ( the area of the usual ball with the radius $\epsilon$ in a plane). In this case,  we denote
$B_{M}(x,r)$ by $B'(x, \epsilon)$. We use the symbol $B'(\epsilon)$, when the center $x$ is not specified . \\\\
The following simple corollary about number of the balls covering
a curve or a surface will be useful.\\\\
{\bf Corollary 2.3.} \\{\it (1) Let  $l$ be the length of a differentiable plan curve $\gamma$ and let $\epsilon>0$.
  There is  a covering of $\gamma$ with a collection  $\{ B_{1}(\epsilon),...,B_{N}(\epsilon)\}$  of balls
     centered on $\gamma$, such that then $N\leq \frac{l}{ \epsilon}$.\\
     (2) Let $M$ be a compact surface in $R^{3}$ with
the area $vol(M)$. There is a covering of $M$ with a collection
$\{ B'_{1}(\epsilon),...,B'_{N}(\epsilon)\}$ of balls of $M$,
such  that $N \leq \frac{3vol(M)}{\pi \epsilon^{2}}$.}
 \begin{proof}
(1): Given a ball $B_{i}(\epsilon)$ centered on $\gamma$. Put
$\gamma_{i}= \gamma \cap B_{i}(\epsilon)$ and let $l_{i}$ be the
length of $\gamma_{i}$. It is clear that $l_{i}$ is at least
equal to
 the diagonal of $B_{i}(\epsilon)$. Then,  $l_{i}\geq 2 \epsilon$.
 We can choose the collection of balls  in such a way that for each $i$, $1<i<N-1$,
  $\gamma_{i}$  has overlaps at most by $\gamma_{i-1}$
 and $\gamma_{i+1}$. Then,
$\sum_{i}l_{i}\leq 2l$ and
\[2l\geq  \sum_{i} l_{i} \geq 2 N \epsilon .\]
(2):  It is possible to choose a collection $\{
B'_{1}(\epsilon),...,B'_{N}(\epsilon)\}$ of  balls of $M$ covering
$M$  in such a way that each point $x \in M$ belongs to at most
three balls. Then, $\sum_{i}vol(B'_{i}(\epsilon)) \leq 3 vol(M)$.
Thus,
\[ N\pi \epsilon^{2} \leq 3vol(M).\]
which gives the result.
\end{proof}
{\bf Lemma 2.4.} {\it Let $\gamma_{1}, ..., \gamma_{n}$ be
sub-curves of $\gamma$, $l$ be the length of $\gamma$ and $l_{i}$
be the length of $\gamma_{i}$ and suppose that for a constant
positive number $a$, $l_{i}>a$.  If $B=\{x_{1},...,x_{m}\} $ be a set of   i.i.d. uniformly distributed points of  
$\gamma$, $0<p\leq 1$, and $m>
\frac{1}{2}(1+\sqrt{1+\frac{8log(1-p)-logn}{log(1-(\frac{a}{l})^{2})}})$,
then with the probability at least  $p$, for all $i \in
\{1,...,n\}$, $\gamma _{i}$ contains at least two different
points of $B$.}
\begin{proof}
Denote by $E_{i}$ the event that no couples of points of $B$
belongs to $\gamma_{i}$. Consider two points $y_{1}, y_{2}$ in
$B$. The probability that  $y_{j}$, $j=1,2$, belong to
$\gamma_{i}$ is $\frac{l_{i}}{l}$ and the probability that both
of the points $y_{1}, y_{2}$ belong to $\gamma_{i}$ is
$(\frac{l_{i}}{l})^{2}$. The number of the different couples
$\{y_{1}, y_{2}\}$ in $B$ is equal to $\binom{m}{2}$. Thus, with
the following probability, no couples of $B$ belong to
$\gamma_{i}$
\[ P(E_{i})=(1-(\frac{l_{i}}{l})^{2})^{\binom{m}{2}}<(1-(\frac{a}{l})^{2})^{\binom{m}{2}}.\]
The probability that for some $i$, $\gamma_{i}$ contains no
couples of $B$ is equal to $P(\bigcup_{i}E_{i})$. What we are
looking is $1-P(\bigcup_{i}E_{i})$. We have
\[ 1-P(\bigcup_{i}E_{i})\geq 1-(P(E_{1})+P(E_{2})+...+P(E_{n})) \geq 1-n(1-(\frac{a}{l})^{2})^{\binom{m}{2}}.\]
A simple computation shows that
$1-n(1-(\frac{a}{l})^{2})^{\binom{m}{2}}>p$ is equivalent to
$\binom{m}{2}>\frac{log(1-p)-logn}{log(1-(\frac{a}{l})^{2})}$. In
other way, if $b$ is a positive number then $\binom{m}{2}\geq b$
iff $m \in (-\infty, \frac{1-\sqrt{1+8b}}{2}] \cup [
\frac{1+\sqrt{1+8b}}{2}, \infty)$, which for positive $m$,
implies $m\geq  \frac{1+\sqrt{1+8b}}{2}$. Putting
$b=\frac{log(1-p)-logn}{log(1-(\frac{a}{l})^{2})}$ gives the
result.
\end{proof}

{\bf Remark 2.5.} The criterion of the previous lemma, by a similar proof, can be mentioned in a general form as follows:\\
Let $A_{1},...,A_{n}$ be subsets of a measurable set $X$ and $B=\{x_{1},...,x_{m}\} \subset X$. If the probability that a point of $X$ belongs to $A_{i}$ is
at least $\alpha$ and $m> \frac{1}{2}(1+\sqrt{1+\frac{8log(1-p)-logn}{log(1-\alpha^{2})}})$, then with the probability at least $p$, each $A_{i}$ contains two different points of $B$.\\\\

{\bf Definition 2.6.} Let $\gamma$ be a connected curve and $x,
x', x'' $ be three different points on $\gamma$.  We say that $x$
is between $x', x''$  and we denote by $[x', x, x'']$, if  for   a
parametrization $\gamma(t)$, $t \in I \subset R$,   of $\gamma$,
\[ x'=\gamma(t'), x=\gamma(t), x''=\gamma(t''),  \ \ t'<t<t'' \ \ or \ \ t''<t<t'.\]

Suppose that we have a point cloud $\{x_{1},...,x_{m}\}$ of the
points on
 a curve $\gamma$, drawn in i.i.d fashion.  The following theorem proposes
a method to estimate(learn) the curvature, with a suitable
confidence by
 choosing sufficiently big number of points. If $x ,y$ are two points on $\gamma$, we denote
 by $d_{\gamma}(x,y)$ the length of the $\gamma$ between $x$ and $y$, which is by definition, the length
 of the shortest path on $\gamma$ joining $x$ to $y$. It is clear
 that $d(x,y)\leq d_{\gamma}(x,y)$.
  \\\\
  {\bf Theorem 2.7.} {\it Let $\epsilon>0$, and  $\gamma$ be a
plan curve with the length $l$  such that $l\geq \epsilon$.  If we choose a set
$B=\{x_{1},..,x_{m}\}$   of   i.i.d. uniformly distributed points of  
$\gamma$ in such a way that \[m>
\frac{1}{2}(1+\sqrt{1+\frac{8log(1-p)-logl+log2\epsilon}{log(1-\frac{\epsilon^{2}}{l^{2}})}}),\]
then,
 with the probability at least  $p$, for each point $x \in \gamma$, there exists $x', x''$ in $B$
 such that $[x', x, x'']$ and $max\{ d(x,x'), d(x,x'')\}<\epsilon$. Thus,   $r^{-1}(x,x',x'')$
  is a $\epsilon$-triangle estimation of $\kappa(x)$.}
\begin{proof}

Given $x \in \gamma$, consider a collection
$A_{i}=B_{\gamma}(y_{i}, \frac{\epsilon}{2})\}$, $1\leq i \leq n$,
of $\gamma$-balls covering $\gamma$ such that
 $x=y_{1}$. Without lose of generality we can assume that for all $i$, $l_{i}=l(A_{i})= \epsilon$ (note that  if $\gamma$ is open bonded and the distance between  $x$ and one of the endpoints, say $e$, is equal to  a number $\delta$ less than  $\frac{\epsilon}{2}$,  we can replace  $\gamma$ by its extension to a longer curve $\bar{\gamma}$ from the endpoint $e$, in such a way that
 $l(\bar{\gamma})=l(\gamma)+(\epsilon-\delta$). Then, assumption $l(A_{i})=\epsilon$ makes no restrictions on our computations). Thus, the probability
 that a chosen point of $\gamma$  belong to $A_{i}$ is  equal to $\alpha=\frac{\epsilon}{l}$.  By Remark 2.5,
if \[m>
\frac{1}{2}(1+\sqrt{1+\frac{8log(1-p)-logn}{log(1-\frac{\epsilon^{2}}{l^{2}})}})
\ \ (*),\] then with the probability at least $p$, there are two
points of $\{x_{1},...,x_{m}\}$ say $x_{\alpha}, x_{\beta}$
belonging to $A_{1}$. If $[x_{\alpha},x, x_{\beta}]$ we are done,
if not then $[x, x_{\alpha}, x_{\beta}]$ or  $[x,
x_{\beta},x_{\alpha}]$. Let $[x,x_{\alpha},x_{\beta}]$ ( the
other case is similar) , consider the
point $z \in \gamma$ such that $[z, x, x_{\alpha}]$ and $d_{\gamma}(z,x)=\epsilon$.\\
By similar arguments, there are two points $x''',x''''$ in
$\{x_{1},...,x_{m}\}$ belonging to $B_{\gamma}(z,
\frac{\epsilon}{2})$.
It is clear that $x''',x'''' \in B_{\gamma}(x,\epsilon)$ and $[x''', x, x_{\alpha}]$.\\
Since $\{A_{i}: 1\leq i\leq n\}$ covers $\gamma$, then  $n\geq
\frac{l}{2 \epsilon}$ which by (*), yields to our result.
\end{proof}

{\bf Remark 2.8.} If $l\geq \epsilon$, then the quantity under the square root in the previous theorem is positive. Because,  by a simple computation we see that 
$1+\frac{8log(1-p)-logl+log2\epsilon}{log(1-\frac{\epsilon^{2}}{l^{2}})}\geq 0$ is equivalent to $l^{3}-(1-p)^{8}2 \epsilon l^{2}+(1-p)^{8}2 \epsilon^{3} \geq 0$.
Put $a=(1-p)^{8}2 \epsilon$, $b=(1-p)^{8}2 \epsilon^{3}$ and consider the function $f(x)=x^{3}-ax^{2}+b$. We show that $f(l)\geq 0$. Since  $f'(x)=3x^{2}-2ax$ then for $x\geq \frac{2a}{3}$, $f'(x)$ is positive. Thus, for $x \geq \frac{2a}{3}$, we have $f(x)\geq f(\frac{2a}{3})$. In other way, after some simple computations, we get
\[ f(\frac{2a}{3})=f(\frac{2(1-p)^{8}2 \epsilon}{3})=2\frac{\epsilon^{3}}{27}(1-p)^{8}(-16(1-p)^{16} +27).\]
Since $0<p\leq 1$ then  $f(\frac{2a}{3})\geq 0$. Now, notice that $l\geq \frac{2a}{3}$ is equivalent to $l \geq \frac{4(1-p)^{8}}{3} \epsilon$, which ( by assumption $l\geq \epsilon$) is true.  Thus, $f(l)\geq  f(\frac{2a}{3}) \geq0$.\\\\
Let $e=(x_{0},y_{0})$ be a unit vector at the origin of $R^{2}$.
The following subset of $R^{2}$ is called a standard cone with the
 radius $r$, angle $\theta$, and direction  $e$.
 \[c_{0}(e,\theta,r)=\]\[\{t(x_{0}cos\psi-y_{0}sin \psi, x_{0}sin
 \psi+y_{0}cos\psi):   t \in [0,r],  \psi \in [-\frac{\theta}{2},
 \frac{\theta}{2}]\}.\]
 We recall that if $M$ is a surface and $x\in M$, then the exponential map  $exp_{x}:T_{x}M \to M$, defines a
 map from $T_{x}M$ (the vector space all vectors tangent to $M$ at
 $x$) to $M$. For each vector $v \in
 T_{x}M$, if $\gamma_{v}$ is the unit speed geodesic in $M$ starting from $x$, then
 $exp_{x}(v)=\gamma_{v}(1)$.\\\\
{\bf Definition 2.9.} Let $M$ be a surface, $x \in M$ and $e$ be
a unit vector in $T_{x}M (\simeq R^{2})$. Consider the exponential
map $exp_{x}:T_{x}M \to M$. The image of a standard cone under the
exponential map is called a cone with the vertex $x$ in $M$. That
is: \[  C=cone(x,e, \theta,r)=exp_{x}(c_{0}(e, \theta,r)).\]
 We denote the  $cone(x,-e, \theta,r)$ by $-C$.
  We say that two points $x_{1}, x_{2}$ of $M$ are well ordered with
   respect  to $cone(x,e, \theta,r)$ if there is a differentiable
   curve
   $\gamma: I \to C \cup(-C) \subset M$ such that $[x_{1}, x,x_{2}]$ on $\gamma$.\\\\
   {\bf Remark 2.10} (see [11]). Let $M$ be a surface and $x\in M$. Consider a
    vector $N_{x}$ which is normal to $M$ at $x$. Planes containing $N_{x}$
    are called normal planes at $x$. The intersection of a normal plane and $M$ will
    form a curve $\gamma$, called a normal section, and the curvature of this curve is the normal curvature.
     The maximum and minimum values of the normal curvatures, say $\kappa_{1}$, $\kappa_{2}$, are called the
     principal curvatures. The normal sections $\gamma_{1}$ and $\gamma_{2}$  whose curvatures are $\kappa_{1}$ and $\kappa_{2}$ are 
     called the principal curves at $x$.
There are (unit) vectors tangent to $M$ at $x$ such that $e_{1}$
is normal to $e_{2}$ and
$\gamma_{1}$, $\gamma_{2}$ are tangent to $e_{1}$, $e_{2}$ at
$x$. $e_{1}$ and $e_{2}$ are called the principal vectors.
The (Gaussian) curvature of $M$ at $x$ is the product  $\kappa= \kappa_{1}\kappa_{2}$.\\\\
{\bf Definition 2.11.} Let $e_{1}$ and $e_{2}$ be the principal
vectors of $M$ at $x$. The cones $C=con(x,e_{1}, \theta, r)$ and
$D=con(x,e_{2},\theta,r)$ are called the principal
$(\theta,r)$-cones at $x$. We say that four distinct points
$x_{1}, x_{2}, x_{3}, x_{4}$ are $(\theta, r)$-principal ordered
at $x$, if
two points are well ordered with respect to $C$ and the other two points are well ordered with respect to $D$.\\\\
{\bf Theorem 2.12.} {\it Let $\epsilon, \theta >0$ and $M$ be a
surface in $R^{3}$ with the area $s$.  If we choose a subset
$B=\{  x_{1},...,x_{m}\}$  of   points, i.i.d. uniformly distributed on $M$, in such away that
\[ m\geq log(\frac{1-p^{\frac{\pi \epsilon^{2}}{12s}}}{4})(log(1-\frac{\theta
 \epsilon^{2}}{2s}))^{-1},\] then with the probability at least $p$, for each point $b \in M$, there are   points $x_{11},x_{12},x_{21},x_{22}$ in $B$ which are $(\theta, \epsilon)$-principal ordered at $b$.}
\begin{proof}
 Consider a covering of $M$ by the
$\frac{\epsilon}{2}$-balls of $M$, say
 \[\{ B'(b_{1},\frac{ \epsilon}{2}),...,B'(b_{n}, \frac{\epsilon}{2})\}, \ \  b_{1}=b.\]
  By Corollary 2.3, we can choose the balls in
  such a way that \[n\leq \frac{3s}{ \pi (\frac{\epsilon}{2})^{2}}=\frac{12s}{\pi \epsilon^{2}}.\]
   Denote by $C_{i}$ and $D_{i}$ the $(\theta, \frac{ \epsilon}{2})$-principal
   cones at $b_{i}$. The area of $C_{i}$ ( and similarly
   $D_{i}$), which we denote by $\mu(C_{i})$ is at least equal to
   the area of a standard $(\theta, \frac{\epsilon}{2})$-cone of
   $R^{2}$. Thus,
\[ \mu(C_{i})\geq \frac{\theta}{2} \epsilon^{2}.\]
Similarly, $\mu(D_{i})\geq \frac{\theta}{2} \epsilon^{2}$. Thus,
the probability that a point of $M$ belongs to $C_{i}$ is
 \[\frac{\mu(C_{i})}{s} \geq \frac{ \theta \epsilon^{2}}{2s}.\]
Let $E_{i}$ be the event that $C_{i}$ contains no points of $B$.
We have
 \[ P(E_{i})\leq( 1- \frac{ \theta \epsilon^{2}}{2s})^{m}.\]
 Similar argument is valid for $-C_{i}$, $D_{i}$ and $-D_{i}$.\\
 Then, for a fixed index $i$, the probability that at least one of the four sets $C_{i}, -C_{i}, D_{i}, -D_{i}$ contains no point of $B$ is
 \[ at \ \ most \ \  4(1- \frac{ \theta \epsilon^{2}}{2s})^{m}.\]
 Consequently, the event that the mentioned four sets  contain a point of $B$ has the probability
 \[at \ \ least \ \ 1- 4(1- \frac{ \theta \epsilon^{2}}{2s})^{m}.\]
 Therefore,  the probability $P$ of the event that for all $i$, the sets\\ $C_{i}, -C_{i}, D_{i}, -D_{i}$,
  contain a point of $B$ is
 \[ P \geq (1- 4(1- \frac{ \theta \epsilon^{2}}{2s})^{m})^{n}.\]
 Thus, $P\geq p$  if
 $ (1- 4(1- \frac{ \theta \epsilon^{2}}{2s})^{m})^{n}\geq p$,
 which is equivalent to
\[mlog(1-\frac{\theta \epsilon^{2}}{2s})\leq
log(\frac{1-\sqrt[n]{p}}{4}).\]
 Since $\theta$ and $\epsilon$ are small,  then
 \[ m\geq log(\frac{1-\sqrt[n]{p}}{4})(log(1-\frac{\theta
 \epsilon^{2}}{2s}))^{-1}.\]
 Since $n\leq \frac{12s}{\pi \epsilon^{2}}$ and for  all $logA$ in the computations we have $A<1$, then

\[log(\frac{1-p^{\frac{\pi \epsilon^{2}}{12s}}}{4})(log(1-\frac{\theta
 \epsilon^{2}}{2s}))^{-1}\geq log(\frac{1-\sqrt[n]{p}}{4})(log(1-\frac{\theta
 \epsilon^{2}}{2s}))^{-1}.\]

Therefore, $P\geq p$ if
\[ m\geq log(\frac{1-p^{\frac{\pi \epsilon^{2}}{12s}}}{4})(log(1-\frac{\theta
 \epsilon^{2}}{2s}))^{-1}.\]

\end{proof}

{\bf Remark 2.13.} (curvature estimation on a surface). Let $M$
be a surface and $b \in M$, and let $x_{11}, x_{12}, x_{21},
x_{22}$ be $(\theta, \epsilon)$-ordered at $p$ such that
$x_{11}$, $x_{12}$ be well ordered with respect to $con(b,e_{1},
\theta, \epsilon)$ and $x_{21}, x_{22}$ be well ordered with
respect to $con(b, e_{2}, \theta, \epsilon)$. Let $\beta_{1}$ be a
broken-geodesic in $M$ ( probably broken at $b$) joining $x_{11}$ to $b$ and $b$ to $x_{12}$ 
( i.e, $[x_{11},b,x_{12}]$ on $\beta_{1}$). Similarly consider the
broken-geodesic $\beta_{2}$ for $x_{21}, x_{22}$. $\beta_{1}$ and
$\beta_{2}$ are estimations for the principal curves at $b$. 
\[  \kappa_{1}(b)\simeq (r(x_{11},b,x_{12}))^{-1}\]
and
\[ \kappa_{2}(b) \simeq (r(x_{21},b,x_{22}))^{-1}.\]
Therefore,
\[ \kappa(b) \simeq  (r(x_{11},b,x_{12}))^{-1}(r(x_{21},b,x_{22}))^{-1}.\]

 Now, by Theorem 2.12 and Remark 2.13, we get the following theorem.\\\\
{\bf Theorem 2.14.}  {\it Let $\epsilon, \theta >0$ and $M$ be a
surface in $R^{3}$ with the area $s$.  If we choose a subset
$B=\{  x_{1},...,x_{m}\}$  of points,  i.i.d. uniformly distributed  on $M$, in such away that
\[  m\geq log(\frac{1-p^{\frac{\pi \epsilon^{2}}{12s}}}{4})(log(1-\frac{\theta
 \epsilon^{2}}{2s}))^{-1},\] then with the probability at least $p$, for each point $b \in M$, there are   points
$x_{11}, x_{12}, x_{21}, x_{22}$ in $B$ such that

\[ \kappa(b) \simeq  (r(x_{11},b,x_{12}))^{-1}(r(x_{21},b,x_{22}))^{-1}.\]
}

 \section{Algorithms}
 There are many algorithms for estimation of the  length  or area of the
 curves or surfaces from a point cloud, which can be used to estimation of $l$ or $s$ in the following
 algorithms.\\\\
{\bf 3.1. Curvature of a curve at a given point.}\\
 1. Input the point $a$, the positive (small) number $\epsilon$, and the expected probability $p$.\\
2. Input the point cloud $B=\{ x_{1},...,x_{m}\}$.\\
3. Compute the length of the curve $l$. Put $\epsilon_{1}=min\{l, \epsilon\}$.\\
4. Compute $\lambda$:\\
\[\lambda=\frac{1}{2}(1+\sqrt{1+\frac{8log(1-p)-logl+log2\epsilon_{1}}{log(1-\frac{\epsilon_{1}^{2}}{l^{2}})}}).\]
5. If $m>\lambda$, then go to 7.\\
 6. Add more points to $B$, then go to 3.\\
 7. Find $y_{1}, y_{2}$ in $B$ such that $[y_{1}, a,
y_{2}] $ and
$max\{ d(a,y_{1}), d(a,y_{2}) \} < \epsilon_{1}$.\\
8. Compute the $\epsilon_{1}$ triangle estimation of $\kappa(a)=r^{-1}(y_{1},a,y_{2})$ (since $\epsilon \leq \epsilon_{1}$, it is also  $\epsilon$ triangle estimation). \\\\
{\bf 3.2. Principal curvature of a surface at a point.}\\\\
Following Remark 2.13 and Theorem 2.14, denote the point cloud by $C$. For estimation of the principal curvatures, 
first step is to find the best estimate for $\beta_{1}$ and its curvature at the point $a$.
 Consider a point $q$ in $C$, near $a$. Let   $q'$ be a point in $C$ with the property that
 $d(q,a)+d(q',a)-d(q,q')$ is equal to  the minimum of the  following set:
  \[\{d(y_{1},a)+d(y_{2},a)-d(y_{1},y_{2}):  \ \ y_{1}, y_{2} \in C\}.\]
   We call $q'$ a $C$-conjugate of $q$.
  If $\beta$ is a geodesic joining $q$ to $a$ and $a$ to $q'$, then its curvature is almost equal to $r^{-1}(q,a,q')$. 
  Thus, we have
  \[ \kappa_{1}(a)=min\{ r^{-1}(c,a,c'): c \in C\}\]
  
  For estimation of $\kappa_{2}(a)$, note that $\beta_{2}$ is included in a plane normal to the plane which
  contains $\beta_{1}$. Thus, if $q$ belongs  to $\beta_{1}$, a point $b$ belongs to $\beta_{2}$ if the vector $ba$ is parallel  to $qa \times q'a$.  Thus, a point $b$ in the point cloud is almost on $\beta_{2}$ if the following number is maximal
  \[ \frac{|<ba,qa \times q'a>|}{|ba||qa \times q'a|}, \ \ b \in C \ \ (*).\]
  Thus,  the second principal curvature
  at $a$ is estimated as:
  \[ \kappa_{2}(a)=min\{ r^{-1}(b,a,b'): b \in C \ \  with \ \ the \ \ property \ \ (*)\}.\]
  Therefore, we can propose the following algorithm for estimation of the principal curvatures and
  the Gaussian curvature of the surfaces using a point cloud on the surface.\\\\
 
1. Input the point $a$, small number $\epsilon$, and
the
expected probability $p$.\\
2. Input the point cloud $C=\{ x_{1},...,x_{m}\}$.\\
3. Compute the area of the surface $s$.\\
4. Compute
\[ \lambda=  log(\frac{1-p^{\frac{\pi \epsilon^{2}}{12s}}}{4})(log(1-\frac{
 \epsilon^{3}}{20s}))^{-1}.\]
 In Theorem 2.14, for simplicity take $\theta=\frac{\epsilon}{10}$.\\
5. If $m>\lambda$, then go to 7.\\
 6. Add more points to $C$, then go to 3.\\
7. Find the collection $D$ of all points $\{y\}$ of
$C$ such that
$d(y,a) < \epsilon$. \\
For each $y \in D$, put:\\
\[  E_{y}=\{d(y,a)+d(c,a)-d(y,c): c \in C\}\]
Put 
\[ E'_{y}=\{ c \in C: d(y,a)+d(c,a)-d(y,c)=min(E_{y})\}\]
\[ E''_{y}=\{ r^{-1}(y,a,c): c \in E'_{y}\}\]
\[k_{y}=min(E''_{y})\]
The first principal curvature is:
\[\kappa_{1}(a)=min\{k_{y}: y \in D\}.\]
8. Select a point $q \in C$ with the property that 
\[ r^{-1}(q,a,q')=\kappa_{1}(a).\]
Put 
\[ F=\{ \frac{|<ba,qa \times   q'a>|}{|ba||qa \times q'a|}: b \in C , b\neq a, qa\neq q'a \neq 0\}.\]
\[ F'=\{ b \in C: \frac{|<ba,qa \times   q'a>|}{|ba||qa \times q'a|}=max{F}\}.\]
\[F''=\{ r^{-1}(b,a,b'): b \in F'\}.\]
The second principal curvature is:
\[ \kappa_{2}(a)=min(F'').\]

9. Return 
  Gaussian curvature, and  Mean curvature at $a$ respectively, as:

 \[  \kappa_{1}(a)\kappa_{2}(a), \ \ \frac{\kappa_{1}(a)+\kappa_{2}(a)}{2}.\]\\\\
{\bf Remark 3.1.} In  estimating   curvature  using a finite set of sample  points, in general, it is not possible to establish an upper bound for the  estimation error.
In fact, for each finite set $B=\{ x_{1},...,x_{n}\} $ of the points in $R^ {m}$  and each positive number $q$, there is a submanifold  $M$ of $R^{m}$,
such that $B \subset M$ and the sectional curvature of $M$ at some point is bigger than $q$r. This is a weak point for all estimations based on the selected finite  set of the  points on manifolds. Nevertheless, our algorithms appear to be effective in practical applications.\\\\
{\bf Remark 3.2.} All necessary codes were prepared by the author( using Python). As shown in the following tables,  the numerical results confirm that the algorithms provide accurate and reliable estimations across various test cases. With the implemented codes, a finite set of points  can be randomly selected and utilized, in conjunction with our
proposed algorithms, to estimate the curvatures. In our experiments, $\epsilon$ and $p$ were set to $0.1$. Reducing the parameter $\epsilon$ can lead to smaller estimation errors.\\
{\tiny
\begin{table}[ht]
\caption{Estimation of the curvature of some curves with  $\epsilon=0.1$ and $p=0.1$} 
\centering 
\begin{tabular}{c c c c} 
\hline\hline 
The curve & Point & The curvature at Point & The estimated curvature at Point \\ [0.5ex] 
\hline 
$x^{2}+y^{2}=25$ & (4,3) & 0.200000  & 0.200008 \\ 
$y=x^{3}+2x$ & (1,3) & 0.004333 & 0.004348 \\
$y=5x^{4}+2x$ & (1,7) & 0.005617 & 0.005762 \\
$x^{2}+y^{2}=0.0025$ & (0,0.05) & 20.000000 & 19.999999 \\
$y=ln(x^{4}+1)+5x^{2}$ & (0,5) & 10.000000 & 9.740785\\ [1ex] 
\hline 
\end{tabular}
\label{table:nonlin} 
\end{table}}
{\tiny
\begin{table}[ht]
\caption{Estimation of the Gaussian curvature of some surfaces with  $\epsilon=0.1$ and $p=0.1$} 
\centering 
\begin{tabular}{c c c c} 
\hline\hline 
The surface & Point & The curvature at Point & The estimated curvature at Point \\ [0.5ex] 
\hline 
$x^{2}+y^{2}+z^{2}=25$ & (0,0,5) & 0.040000  & 0.039999 \\ 
$z=x^{3}+2x+y^{2}+y$ & (1,2,9) & 0.004618& 0.004204 \\
$z=\sqrt{x^{2}+y^{2}}$ & $(1,2,\sqrt{5})$ & 0.027777 & 0.025705 \\
$z=x^{2}+0.25y^{2}$ & (1,2,2) & 0.027777 & 0.0240704 \\
 \\ [1ex] 
\hline 
\end{tabular}
\label{table:nonlin} 
\end{table}}
------------------

\end{document}